\documentclass[12pt]{article}
\usepackage{amsmath}
\usepackage{amssymb}
\usepackage{graphicx}
\usepackage{tikz}
\usepackage{pdfsync}

\newcommand{\mat}[1]{\ensuremath{
\left[\begin{matrix}#1
\end{matrix}\right]
}}

\newcommand{\nwline}[1]{\put(#1){\line(-1,1){.5}}}
\newcommand{\swline}[1]{\put(#1){\line(-1,-1){.5}}}
\newcommand{\neline}[1]{\put(#1){\line(1,1){.5}}}
\newcommand{\seline}[1]{\put(#1){\line(1,-1){.5}}}
\newcommand{\veeup}[1]{\nwline{#1}\neline{#1}}
\newcommand{\veedown}[1]{\swline{#1}\seline{#1}}
\newcommand{\e}{\varepsilon}
\newcommand{\sgn}{{\rm sgn}}
\newcommand{\Hom}{{\rm Hom}}
\newcommand{\Ext}{{\rm Ext}}

\newcommand{\spot}[1]{\put(#1){\put(-.1,-.1){$\bullet$}}}

\begin{document}

\pagestyle{empty}

\centerline{\bf Mixed Cobinary Trees}\bigskip

\centerline{Kiyoshi Igusa and Jonah Ostroff}


\begin{abstract}
We develop basic cluster theory from an elementary point of view using  a variation of binary trees which we call mixed cobinary trees. We show that the number of isomorphism classes of such trees is given by the Catalan number $C_n$ where $n$ is the number of internal nodes. We also consider the corresponding quiver $Q_\epsilon$ of type $A_{n-1}$. As a special case of more general known results about the relation between $c$-vectors, representations of quivers and their semi-invariants, we explain the bijection between mixed cobinary trees and the vertices of the generalized associahedron corresponding to the quiver $Q_\epsilon$.
\end{abstract}
\bigskip

\centerline{\bf Introduction}
\smallskip

In a sequence of very influential papers, Fomin and Zelevinsky developed the theory of cluster algebras. Here we look at a concept introduced in [FZ], namely the \emph{$c$-vectors} which form the columns of the bottom half of the exchange matrix $\tilde B$. (See \eqref{eq: def of exchange matrix} below.) Instead of starting with the general definition of cluster mutation, we will examine the elementary concept of \emph{mixed cobinary trees} and study how these trees change when the slope of an internal edge changes. We develop properties of these trees with the aim of explaining known relations between different cluster notions, in particular the relation between representations of quivers, their semi-invariants and corresponding $c$-vectors.

This paper seeks to explain cluster concepts in elementary terms. We start with the definition of {mixed cobinary trees}. We show in Theorem 1 that there are a Catalan number $C_n$ of these trees (up to isomorphism) with $n$ internal nodes. Then we construct a partition of Euclidean space $\mathbb R^n$ with each part $R(T)$ supporting one tree. \emph{Mutation} of trees is defined to be the process of moving from one part to another adjacent part through a well-defined wall. The combinatorics of tree mutation is studied and encoded in the exchange matrix $\tilde B(T)$ using a formula of Nakanishi and Zelevinsky [NZ]. Theorem 2 shows that $\tilde B(T)$ mutates according to the Fomin-Zelevinsky definition of cluster mutation which we review.

A \emph{quiver} is a directed graph and a \emph{representation} of a quiver is given by assigning a vector space to each vertex and a linear map to each arrow. We briefly describe the notion of \emph{clusters} of representations (introduced in [BMRRT]) using elementary matrix equations. In Theorem 3, we restate the virtual stability theorem of [IOTW]. Roughly, this says that, up to sign, the classical $c$-vectors associated to a cluster are linearly related to the \emph{weights of virtual semi-invariants} associated to the cluster. Finally, Theorem 4 describes the bijection between clusters for $Q_\e$ and MCT's and implies (Corollary 2) that the bijection preserves $c$-vectors, i.e., the classical $c$-vectors of a cluster are equal to the $c$-vectors of the corresponding tree. We give two example to explain these statements. Example 1 illustrates the definition of virtual semi-invariants, and the stability conditions. Example 2 illustrates the statement of Theorem 4.

Theorem 4 could be proved using Theorem 2 since the $c$-vectors are uniquely determined by the fact that they satisfy the Fomin-Zelevinsky mutation rules. However, we give a more direct proof of this correspondence using Theorem 3. Our proof of Theorem 4 is based on Theorem 1. We show that there is an epimorphism from the set of clusters of representations to the set of mixed cobinary trees. Using the well-known fact that there are a Catalan number of clusters, we conclude that this mapping is a bijection. The epimorphism is given by showing that the boundaries of the regions $R(T)$ satisfy the stability conditions and therefore form a subset of the walls separating the regions corresponding to clusters. There is also a more technical summary at the end of the paper. 

In a subsequent paper, with different authors, these results will be extended to the affine $\tilde A_{n-1}$ case.


\bigskip

\centerline{\bf Counting mixed cobinary trees}\smallskip

A binary tree (sometimes called a ``full binary tree'') is a rooted tree in which every internal vertex has exactly two children, and whose vertices are endowed with two orderings: a left-to-right total ordering, and a bottom-to-top partial ordering. Both of these orderings are induced by the edges of the tree: a node is to the left of its right descendants, to the right of its left descendants, below its parent, and above its children. This may seem like more information than we need, but it will be useful for the following generalization.

It is well known that the number of binary trees with $n$ internal nodes is the $n$th Catalan number $C_n = \displaystyle\frac{1}{n+1}{2n \choose n}$. For instance, here are the five binary trees with three internal nodes.

%
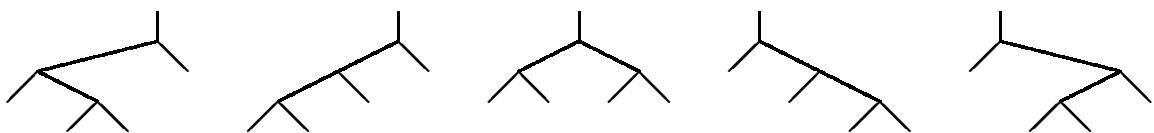
\begin{figure}[htbp]
\begin{center}
%
{
\setlength{\unitlength}{.8cm}
{\mbox{
\begin{picture}(18,2)
      \thicklines
   \put(0,0){
  \veedown{1,.5}
  \swline{0,1}
  \seline{2,1.5}
  \put(2,1.5){\line(0,1){.5}}
  \qbezier(0,1)(.5,.75)(1,.5)
  \qbezier(0,1)(1,1.25)(2,1.5)
   }
\put(4,0){
\veedown{0,.5}
\qbezier(0,.5)(1,1)(2,1.5)
\seline{1,1}
\seline{2,1.5}
\put(2,1.5){\line(0,1){.5}}
   \put(3,0){
   \veedown{1,1}
   \veedown{3,1}
   \qbezier(1,1)(1.5,1.25)(2,1.5)
   \qbezier(2,1.5)(2.5,1.25)(3,1)
   \put(2,1.5){\line(0,1){.5}}
   } 
   \put(8,0){
   \veedown{2,.5}
   \swline{1,1}
   \swline{0,1.5}
   \put(0,1.5){\line(0,1){.5}}
   \qbezier(0,1.5)(1,1)(2,0.5)
   } 
   \put(12,0){
   \veedown{1,.5}
   \swline{0,1.5}
   \seline{2,1}
   \qbezier(1,.5)(1.5,.75)(2,1)
   \qbezier(0,1.5)(1,1.25)(2,1)
   \put(0,1.5){\line(0,1){.5}}
   } 
   } 
\end{picture}}
}}
\caption{These are the $C_3=5$ binary trees with 3 internal nodes.}
\label{figureone}
\end{center}
\vspace{-12pt}
\end{figure}

Binary trees arise naturally when studying associative algebras, where each node represents a multiplication of the factors corresponding to its children. To study coalgebras, however, we'll need a new type of tree with a new type of node to represent comultiplication.

Specifically, a {\em mixed cobinary tree} is a tree with two types of internal nodes: $\Lambda$-nodes to represent multiplication, which have one parent (either to the left or right), one left child, and one right child; and V-nodes to represent comultiplication, which have one child (either to the left or right), one left parent, and one right parent. As with regular binary trees, the vertices of a mixed cobinary tree are endowed with two orderings: a left-to-right total ordering, and a bottom-to-top partial ordering. These orderings also respect the structure of the tree: a $\Lambda$-node is to the right of its left descendants, to the left of its right descendants, below its parent, and above its children; a V-node is to the right of its left ancestors, to the left of its right ancestors, below its parents, and above its single child. Here by the ``left descendants'' (and analogously for right descendants and left and right ancestors) of a $\Lambda$-node $\lambda$, we mean all nodes $u$ for which the path from $\lambda$ to $u$ travels through $\lambda$'s left child; it is not the case that all such nodes are necessarily below $\lambda$, as the following example illustrates:

%
\begin{figure}[htbp]
\begin{center}
%
{
\setlength{\unitlength}{1cm}
{\mbox{
\begin{picture}(3.5,2)
      \thicklines
\put(-.5,.5){\line(1,1){.5}}
\put(0,.5){\qbezier(0,.5)(1,.25)(2,0)}
\put(0,.5){\qbezier(0,.5)(.5,.75)(1,1)}
\put(2,0.5){\line(0,-1){.5}}
\put(2,0){\qbezier(0,.5)(.5,.75)(1,1)}
\seline{3,1}
\put(3,1){\line(0,1){.5}}
\put(.5,2){\line(1,-1){.5}\put(0,-.5){\line(1,1){.5}}}
\put(3.1,1){$\lambda$}
\put(-.3,1.2){$t_1$}
\put(1.1,1.3){$t_2$}
\put(2.1,.3){$t_3$}
\end{picture}}
}}
\caption{Nodes $t_1,t_2,t_3$ are left descendants of $\lambda$ but $t_2$ is above $\lambda$.}
\label{figuretwo}
\end{center}
\vspace{-12pt}
\end{figure}
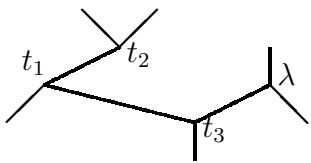
A more technical definition can be given:

{\bf Definition.} A \emph{mixed cobinary tree} is a finite tree $T$ which is linearly embedded in the $xy$-plane satisfying the following conditions.

(1) There are no horizontal edges in $T$.

(2) The internal nodes of the tree have distinct $x$ coordinates which we take to be $1,2,3,\cdots,n$.

(3) The internal nodes are either of \emph{positive} type (also called $\Lambda$) with two children and one parent or of \emph{negative} type (also called V) with one child and two parents.

(4) If $(x,y)$ is a node of $T$ of positive type then the vertical wall $\{(x,z)\,:\,z<y\}$ is disjoint from $T$ and separates the two children of the node.

(5)  If $(x,y)$ is a node of $T$ of negative type then the vertical wall $\{(x,z)\,:\,z>y\}$ is disjoint from $T$ and separates the two parents of the node.

We define two mixed cobinary trees to be \emph{isomorphic} if their underlying tree structure is isomorphic with corresponding inner vertices having identical left-to-right total ordering and bottom-to-top partial ordering. Note that the internal edges of the tree form the Hasse diagram of the bottom-to-top partial ordering. However, unlike regular binary trees is the fact that the above rules do {\em not} induce a unique left-to-right ordering of the vertices of the tree. Consider the following trees, whose local structures are isomorphic but which are endowed with two distinct left-to-right orderings of the vertices.

%
\begin{figure}[htbp]
\begin{center}
%
{
\setlength{\unitlength}{1cm}
{\mbox{
\begin{picture}(7,2)
      \thicklines
\veedown{1,.5}
\veeup{2,1.5}
\swline{0,1}
\qbezier(0,1)(.5,.75)(1,.5)
\qbezier(0,1)(1,1.25)(2,1.5)
\put(-.3,1.2){$t_1$}
\put(2.05,1.2){$t_3$}
\put(1.1,.5){$t_2$}
\put(4.5,0){
\veedown{2,.5}
\veeup{1,1.5}
\swline{0,1}
\qbezier(0,1)(1,.75)(2,.5)
\qbezier(0,1)(.5,1.25)(1,1.5)
\put(-.3,1.2){$t_1$}
\put(1.05,1.2){$t_2$}
\put(2.1,.5){$t_3$}
} 
\end{picture}}
}}
\caption{These trees have distinct left-to-right order of nodes.}
\label{figurethree}
\end{center}
\vspace{-12pt}
\end{figure}
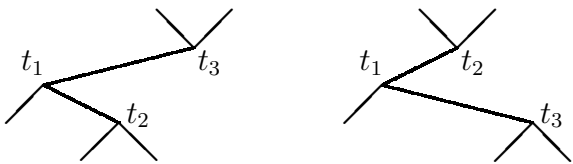

It may seem unnatural, then, to bother endowing a mixed cobinary tree with a total left-right ordering. But notice that in each case where multiple distinct orderings arise from the same tree structure, those orderings in turn yield distinct left-right orderings of $\Lambda$- and V-nodes. For example, in Figure 3, the two node orderings are $(\Lambda,\Lambda,V)$ and $(\Lambda,V,\Lambda)$.

Define the \emph{$\varepsilon$ vector} of a mixed cobinary tree $T$ with $n$ internal vertices to be the $n$-tuple $\varepsilon(T) = (\varepsilon_1,\dots,\varepsilon_n)$, where $\varepsilon_i = 1$ if the $i$-th internal vertex of $T$ from left-to-right is a $\Lambda$-node, and $-1$ if it's a V-node. It turns out that the number of mixed cobinary trees with a particular $\varepsilon$ vector is $C_n$. For instance, here are the 5 MCTs with $\varepsilon$ vector $(-1,-1,1)$.

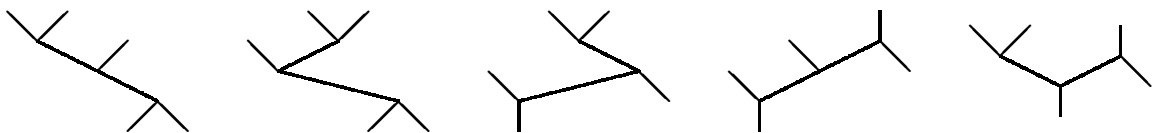
\begin{figure}[htbp]
\begin{center}
%
{
\setlength{\unitlength}{.8cm}
{\mbox{
\begin{picture}(18,2)
      \thicklines
\put(12, 0){ 
\put(0,0){\line(0,1){.5}} 
\nwline{0,.5}
\nwline{1,1}
\seline{2,1.5}
\put(2,1.5){\line(0,1){.5}}
\qbezier(0,.5)(1,1)(2,1.5)
} 
\put(8, 0){ 
\put(0, 0){\line(0,1){.5}}
\nwline{0,.5}
\veeup{1,1.5}
\seline{2,1}
\qbezier(0,.5)(1,.75)(2,1)
\qbezier(1,1.5)(1.5,1.25)(2,1)
} 
\put(4, 0){ 
\veedown{2,.5}
\veeup{1,1.5}
\nwline{0,1}
\qbezier(0,1)(.5,1.25)(1,1.5)
\qbezier(0,1)(1,.75)(2,.5)
} 
\put(0, 0){ 
\veeup{0,1.5}
\veedown{2,.5}
\neline{1,1}
\qbezier(0,1.5)(1,1)(2,.5)
}
\put(16,-.25){ 
\veeup{0,1.5}
\put(2,1.5){\line(0,1){.5}}
\seline{2,1.5}
\put(1,.5){\line(0,1){.5}}
\qbezier(0,1.5)(.5,1.25)(1,1)
\qbezier(2,1.5)(1.5,1.25)(1,1)
} 
\end{picture}}
}}
\caption{There are $C_3=5$ mixed cobinary trees with $\e$ vector $(-1,-1,1)$.}
\label{figurefour}
\end{center}
\vspace{-12pt}
\end{figure}
\smallskip

{\bf Theorem 1:} \emph{Let $\varepsilon \in \{\pm1\}^n$. Then the number of mixed cobinary trees $T$ such that $\varepsilon(T) = \varepsilon$ is the $n$th Catalan number $C_n$.}\smallskip

{\bf Proof:} We'll prove this first by induction, and then by a relatively easy to describe bijection. That the latter is indeed invertible, however, will only be obvious after an observation in the first proof.

The base case is clear: there is $C_1=1$ tree with $\varepsilon(T) = (1)$ and one tree with $\varepsilon(T) = (-1)$. The second tree looks like the letter {\sf Y} and the first tree is the same upside-down: $\ ^{\scalebox{1}[-1]{\text{\sf Y}}}$.

For the inductive step, first note that we need only consider the case where $\varepsilon = (\varepsilon_1,\dots,\varepsilon_{n-1},1)$, since flipping any mixed cobinary tree $T$ over a horizontal axis gives a new cobinary tree $T'$ with $\varepsilon(T')=-\varepsilon(T)$, and by symmetry of the conditions this is a bijection between mixed cobinary trees with those $\varepsilon$ vectors.

Suppose then that $T$ is a tree whose rightmost internal vertex is a $\Lambda$-node $\lambda$, and consider the two (possibly empty) trees obtained by deleting $\lambda$. One of those trees, the ``upper tree'' $U$, reached by traveling from $\lambda$ towards its parent, will have $i$ vertices; the other, the ``lower tree'' $D$, reached by traveling from $\lambda$ towards its left child, will have $n-i-1$ vertices. There are no internal vertices in the direction of $\lambda$'s right child, because $\lambda$ is the rightmost vertex of $T$.

Certainly both of these sub-trees are mixed cobinary trees. We will show that knowing $\varepsilon$, $|U|=i$ and $|D|=n-i-1$ determines $\varepsilon(U)$ and $\varepsilon(D)$, and that each choice of $U$ and $D$ with those $\varepsilon$ vectors yields a unique tree $T$. And since by induction there are $C_i$ choices for $U$ and $C_{n-i-1}$ choices for $D$, the number of choices for $T$ is $$\sum_{i=0}^{n-1} C_i C_{n-i-1} = C_n,$$a well known recurrence for the Catalan numbers.

{\bf Claim:} If $U$ contains a $\Lambda$-node $\lambda'$, then every vertex in $D$ is to the right of $\lambda'$.

{\bf Proof:} There are two cases. In both cases, showing that $\lambda$ is a descendant of a vertex in $U$ implies that it is a right descendant, since $\lambda$ is not to the left of any other vertex. Consider the path from $\lambda'$ to $\lambda$ in $T$. There are two cases: if this path begins with a down step, then $\lambda$ is a descendant of $\lambda'$, so it's a right descendant of $\lambda'$, and therefore so are all vertices in $D$, which means they are to the right of $\lambda'$. If the path begins with an up step, then at some point it must switch from an up step to a down step. This can only happen at another $\Lambda$ node $\lambda''$, so both $\lambda'$ and $\lambda$ are descendants of $\lambda''$. This can only happen if $\lambda'$ is to the left of $\lambda''$, which in turn is to the left of all vertices in $D$. This completes the proof.

{\bf Claim:} If $D$ contains a V-node $v'$, then every vertex in $U$ is to the right of $v'$.

{\bf Proof:} The proof is identical.

{\bf Claim:} If $U$ contains a $\Lambda$-node, then $D$ does not contain a V-node.

{\bf Proof:} Suppose to the contrary that $U$ contains a $\Lambda$-node $\lambda'$ and $D$ contains a V-node $v'$. Then $v'$ is to the right of $\lambda'$ and $\lambda'$ is to the right of $v'$.

{\bf Claim:} Given $\varepsilon$, $|U|$, and $|D|$, the vectors $\varepsilon(U)$ and $\varepsilon(D)$ are uniquely determined.

{\bf Proof:} Given $\varepsilon = (\varepsilon_1,\dots,\varepsilon_{n-1},1)$, let $t_i$ be the $i$-th vertex from the left in $T$ (so $t_i$ is a $\Lambda$ node iff $\varepsilon_i=1$). If $\varepsilon_i = \varepsilon_j = 1$ with $i < j$, then $t_j \in U$ implies $t_i \in U$, or else $t_i$ is a vertex of $D$ to the left of a $\Lambda$-node $t_j$ in $U$, contradicting the previous claim. Likewise if $\varepsilon_i = \varepsilon_j = -1$ with $i < j$, then $t_i \in U$ implies $t_j \in U$, and if $\varepsilon_i = 1$ and $\varepsilon_j = -1$ (for any ordering of $i$ and $j$, then $t_i \in U$ implies $t_j \in U$. These implications give a complete ordering of the possible vertices of $U$, so knowing the number of vertices in $U$ completely determines which vertices they are.

Finally, for any choices of $U$ and $D$ with the appropriate $\varepsilon$ vectors, there is a unique way to arrange their vertices from left to right (given by the above characterization), and a unique way to attach them to $\lambda$ in $T$: $\lambda$ is attached to the rightmost down-leaf of $U$ and the rightmost up-leaf of $D$. It's easy to see that any such pairing of $U$ and $D$ will yield a valid choice for $T$, completing the proof of Theorem 1.\smallskip

{\bf A Bijection:} We can also characterize the correspondence between mixed cobinary trees with a particular $\varepsilon$ vector and regular binary trees through the following bijection: imagine a mixed cobinary tree to be made up of sticks attached loosely attached to each other on the ends as in a child's mobile. Put a tack at the end of the rightmost leaf of the rightmost vertex, then let the entire configuration swing down as affected by gravity, keeping the orientations around each vertex but falling so that the fixed leaf points up and all other leaves point down. The result is an ordinary binary tree $T'$. Applying this bijection to the 5 trees of Figure 4 yields, in order, the 5 trees of Figure 1.

The inversion of this function is not easy to describe, but it does exist: in the case where the rightmost root of $T$ was originally a $\Lambda$-node (the other case is analogous), the left-descendants of the root of $T'$ make up the vertices of $U$, and the right-descendants the vertices of $D$. Knowing the sizes of those subtrees, as before, gives a unique way to shuffle those vertices together from left to right, so it's possible to pick up the tree again and return the nodes to their original orientations.\bigskip

\centerline{\bf Mutation of mixed cobinary trees}\smallskip

If $T$ is a MCT of length $n$ let $\pi(T)$ denote the set of permutations of $n$ given as follows. $\sigma\in S_n$ lies in $\pi(T)$ if the internal nodes of $T$ can be taken to lie at the points $(i,\sigma(i))$ for $i=1,\cdots,n$. For example, in Figure 2 we have $\sigma(4)>\sigma(3)<\sigma(1)<\sigma(2)$. So, $\pi(T)=\{3412, 2413, 2314\}$. \smallskip

{\bf Lemma} (Uniqueness Lemma) \emph{Every permutation of $n$ lies in $\pi(T)$ for a unique MCT $T$.}\smallskip

We prove this by induction on $n$ simultaneously with the following lemma.
\smallskip

{\bf Lemma.} \emph{If a MCT has $p$} V\emph{-nodes then it has $p+1$ parent leaves which are separated by the vertical lines going up from these nodes. Similarly, if there are $q$ $\Lambda$-nodes then there are $q+1$ descending leaves which are separated by the vertical lines descending from the $\Lambda$-nodes.}

%
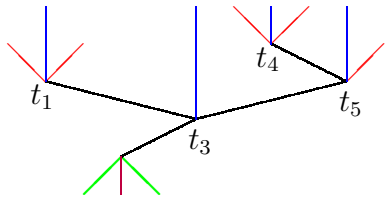
\begin{figure}[htbp]
\begin{center}
%
{
\setlength{\unitlength}{1cm}
{\mbox{
\begin{picture}(6,2.5)
      \thicklines
\color{green}
\veedown{2.1,.5}
}
\put(0,0){
{\color{red}
\veeup{1,1.5}
\veeup{4,2}
\neline{5,1.5}
}
{
\qbezier(1,1.5)(2,1.25)(3,1)
\qbezier(2,.5)(2.5,.75)(3,1)
\qbezier(3,1)(4,1.25)(5,1.5)
\qbezier(4,2)(4.5,1.75)(5,1.5)
\put(2,.5){\color{purple}\line(0,-1){.5}}
{\color{blue}
\put(1,1.5){\line(0,1){1}}
\put(3,1){\line(0,1){1.5}}
\put(4,2){\line(0,1){.5}}
\put(5,1.5){\line(0,1){1}}
}
\put(3.8,1.7){$t_4$}
\put(2.9,.6){$t_3$}
\put(.8,1.2){$t_1$}
\put(4.9,1.1){$t_5$}%
}%
\end{picture}}
}}
\caption{{\color{blue}Vertical ascending walls} starting at V-nodes separate the {\color{red}parent leaves}. {\color{purple}Vertical descending walls} starting at $\Lambda$-nodes separate the {\color{green}children leaves}.}
\label{figureA}
\end{center}
\vspace{-12pt}
\end{figure}

{\bf Proof:} Both lemmas are clearly true for $n=1$. So, suppose $\sigma$ is a permutation of $n\ge2$ and let $k=\sigma^{-1}(n)$. We remove the number $k$ from the set $\{1,\cdots,n\}$. There are two cases.

If $\e_k=1$ then the removal of the $\Lambda$-node $t_k$ for any tree $T$ will break up the tree into two subtrees $T_1,T_2$. These are uniquely given by the permutations $\sigma_1,\sigma_2$ of $k-1,n-k$ respectively given by the total ordering on the sets $\{1,\cdots,k-1\}$ and $\{k+1,\cdots,n\}$ induced by $\sigma$. The tree $T$ is obtained by attaching the rightmost parent leaf of $T_1$ and the leftmost parent leaf of $T_2$ to $t_k$. The new tree $T$ has only one parent leaf between the rightmost ascending wall of $T_1$ and leftmost ascending wall of $T_2$. So, both lemmas hold in this case.

If $\e_k=-1$ then removing the V-node $t_k$ from $T$ will produce a single tree $T'$ given as follows. Take the permutation $\sigma'$ and sign function $\e'$ of $n-1$ given by $\sigma'(i)=\sigma(i)$ and $\e_i'=\e_i$ if $i<k$ and $\sigma'(i)=\sigma(i+1),\e_i'=\e_{i+1}$ if $i\ge k$. By induction on $n$, there is a unique tree $T'$ with $\e(T')=\e'$ with $\sigma'\in \pi(T')$. On this tree take the unique parent leaf which lies between the vertical ascending walls at the nodes $t_i,t_j$ where $i,j$ are maximal and minimal respectively so that $i<k\le j$ and $\e_i'=\e_j'=-1$. Attach the new node $t_k$ to this uniquely determined parent leaf. The new node has a new pair of parents which are separated by the vertical wall ascending from $t_k$. This gives back $T$ satisfying the wall-leaf separation condition. And, $T$ is uniquely determined. So, both lemmas hold.

\smallskip

Given a mixed cobinary tree $T$ we will define a region $R(T)$ in $\mathbb R^n$ by $n-1$ inequalities corresponding to the edges of $T$ and prove that this region is the set of all $x=(x_1,\cdots,x_n)\in\mathbb R^n$ so that the points $(i,x_i)$, for $i=1,\cdots,n$ form an allowable set of internal nodes for $T$. 

As an example, for $T$ as in Figure 2 we define
\[
	R(T)=\{x\in\mathbb R^4\,|\, x_2> x_1> x_3< x_4\}.
\]
In general we take the $n-1$ internal edges $\ell_i$ of $T$ and define $R(T)$ to be given by the inequalities $\sgn (\ell_i)(x_{q_i}-x_{p_i})>0$ where $\sgn (\ell_i)$ is the sign of the slope of the edge $\ell_i$. These inequalities are: $x_{p_i}< x_{q_i}$ for each $\ell_i$ with positive slope and $x_{p_i}> x_{q_i}$ for each $\ell_i$ with negative slope. \smallskip

{\bf Proposition 1:} \emph{For any $x\in\mathbb R^n$, the points $(i,x_i)$, for $i=1,\cdots,n$ form an allowable set of internal nodes for $T$ if and only if $x\in R(T)$.}\smallskip

{\bf Proof:} As the real numbers $x_j$ varies, the node $(j,x_j)$ moves only vertically. By assumption, the slopes of all the internal edges remains correct. So, the only way that a point in $R(T)$ could fail to give a tree isomorphic to $T$ is if the tree is not embedded. 

{\bf Claim:} Take any two point in the tree $T$ with the same $x$-coordinate, say $v,w$, then the unique path from $v$ to $w$ is either always increasing or always decreasing. 

{\bf Proof:} If this path reaches a maximum, say $m$, (or minimum) in the interior of the path then the point $m$ must be a node $(i,x_i)$ with $\e_i=1$. In that case, the rest of the path cannot cross the vertical line descending from $m$. So the endpoints of the path must have distinct $x$ coordinates contradicting our starting assumption. 

As $x$ varies in $R(T)$, this path will always be increasing or decreasing since the signs of the slopes of all internal edges is fixed. So, the tree will always be embedded and the points $(i,x_i)$ form the nodes of a tree isomorphic to $T$. The converse is clear by definition of tree isomorphism. So, Proposition 1 is proved.\smallskip

By the Uniqueness Lemma, Proposition 1 and Theorem 1, we get the following corollary where we use the notation $H$ for the hyperplane in $\mathbb R^n$ of all $x$ so that $\sum x_i=0$. Also, let $\overline R(T)$ be the closure of $R(T)$.\smallskip

{\bf Corollary 1:}
\emph{The intersection $\overline R(T)\cap H$ is the cone of an $n-2$ simplex whose $n-1$ sides are given by the equalities $x_{p_i}= x_{q_i}$. The sets $\overline R(T)$ cover $\mathbb R^n$ and their interiors $R(T)$ are disjoint. Thus, when intersected with the unit sphere $S^{n-2}$ in the hyperplane $H$, we obtain a partition of the sphere into $C_n$ simplices of dimension $n-2$.}
\smallskip

Thus, when the point $x$ moves through one of the walls of the region $\overline R(T)$ given by the equality $x_{p_i}= x_{q_i}$, we enter a new region $R(T^\ast)$ for a new MCT $T^\ast$ which, in Proposition 2 below, we will show is uniquely determined by $T$ and the edge $\ell_i=(p_i,q_i)$ whose slope is being reversed. We say that $T^\ast$ is obtained from $T$ by \emph{mutation in the $i$-direction} and write $T^\ast=\mu_i(T)$. An example is given in Figure 6 below. It is clear that $\mu_i\mu_i(T)=T$ assuming the uniqueness of mutation.

%
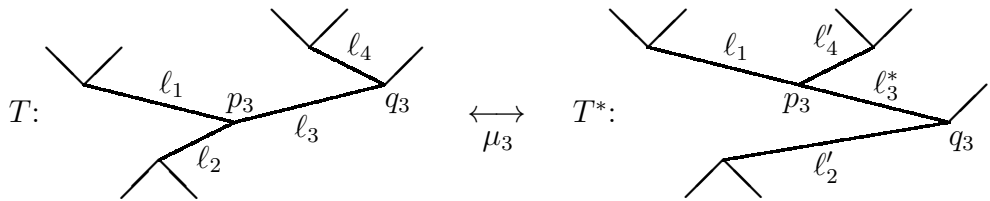
\begin{figure}[htbp]
\begin{center}
%
{
\setlength{\unitlength}{1cm}
{\mbox{
\begin{picture}(12,2.5)
      \thicklines
\veeup{1,1.5}
\veedown{2,0.5}
\veeup{4,2}
\neline{5,1.5}
\qbezier(1,1.5)(2,1.25)(3,1)
\qbezier(2,.5)(2.5,.75)(3,1)
\qbezier(3,1)(4,1.25)(5,1.5)
\qbezier(4,2)(4.5,1.75)(5,1.5)
\put(2,1.4){$\ell_1$}
\put(2.5,.4){$\ell_2$}
\put(2.9,1.2){$p_3$}
\put(3.8,.8){$\ell_3$}
\put(4.5,1.9){$\ell_4$}
\put(5,1.2){$q_3$}
   \put(7.5,0){   \put(0,1){$T^\ast$:}
\veeup{1,2}
\veedown{2,.5}
\veeup{4,2}
\neline{5,1}
\qbezier(2,.5)(3.5,.75)(5,1)
\qbezier(1,2)(2,1.75)(3,1.5)
\qbezier(3,1.5)(3.5,1.75)(4,2)
\qbezier(3,1.5)(4,1.25)(5,1)
\put(2.8,1.2){$p_3$}
\put(5,.7){$q_3$}
\put(2,1.9){$\ell_1$}
\put(3.2,.3){$\ell_2'$}
\put(4,1.4){$\ell_3^\ast$}
\put(3.2,2){$\ell_4'$}
   }
   \put(6.1,1){$\longleftrightarrow$}
   \put(6.3,.7){$\mu_3$}
   \put(0,1){$T$:}
\end{picture}}
}}
\caption{When $\ell_3$, with positive slope becomes $\ell_3^\ast$ with negative slope, the rightmost child of the left endpoint $p_3$ of $\ell_3$ becomes the leftmost child of the right endpoint $q_3$ of $\ell_3^\ast$ and, similarly, the leftmost parent of $q_3$ becomes the rightmost parent of $p_3$.}
\label{figureB}
\end{center}\vspace{-12pt}
\end{figure}
{\bf Remark 1:} For $n\ge3$, the sphere $S^{n-2}$ is path connected. So, it follows from Corollary 1 that any MCT $T$ can be obtained from any other MCT $T'$ by a finite sequence of mutations. Namely, a path in $S^{n-2}$ should be chosen connecting unit vectors in $R(T)$, $R(T')$. By transversality, it will pass through a finite sequence of walls. The corresponding mutations will change $T$ to $T'$.  For $n=1,2$ this statement is also easily seen to be true. We will show in Theorem 4 below that the partition of the sphere given in Corollary 1 agrees with the triangulation of the sphere given in [IOTW], stated below as Theorem 3.
\smallskip

There is a \emph{reverse tree} which is given by taking the mirror image of $T$ through a vertical line. Reversing of trees commutes with mutation and changes the signs of all the edges. So, in order to study the effect of the mutation $\mu_i$, we can assume that $\sgn (\ell_i)=+$. The mirror image of $T$ through a horizontal line is also a mixed cobinary tree with the signs $\e_i$ reversed.\smallskip

{\bf Definition:} We define the \emph{$c$-vector} of $\ell_i$ to be $\sgn (\ell_i)$ times the sum of unit vectors $e_j$ for $p_i\le j<q_i$. We denote this by $c_i(T)$. The \emph{$c$-matrix} of $T$ is defined to be the $n-1\times n-1$ matrix $C(T)$ whose columns are the $c$-vectors $c_i(T)$.\smallskip

For example, if $T$ and $T^\ast=\mu_3(T)$ are as in Figure 6 above, we have
\[
	C(T)=\small\mat{-1&0&0&0\\
	-1&1&0&0\\
	0&0&1&0\\
	0&0&1&-1}
	\quad 
	C(T^\ast)=\mat{-1&0&0&0\\
	-1&1&0&0\\
	0&1&-1&1\\
	0&1&-1&0}
\]
The second matrix is obtained from the first matrix by column operations as follows.

(1) Add Column 3 to Column 2 and to Column 4.

(2) Change the sign of Column 3.

The same instructions in the same order will change the second matrix into the first. The mutation $\mu_k$ adds the $k$th column of $C$ to certain other columns as determined by the following proposition, then changes the sign of the $k$th column of $C$. We note that $T$ is uniquely determined by its $c$-vectors $c_i(T)$ since these determine the edges of $T$ and the signs of their slopes.

\smallskip

{\bf Proposition 2:} \emph{Suppose that $\ell_k$ is an internal edge of $T$ whose left node $p_k$ is below its right node $q_k$, i.e., $\sgn (\ell_k)=+$. Let $\pi\in\pi(T)$ so that $\pi(p_k)+1=\pi(q_k)$. Let $\pi'$ be the composition of $\pi$ with the simple transposition reversing the values of $\pi(p_k),\pi(q_k)$. Let $T^\ast$ be the unique MCT so that $\pi'\in\pi(T^\ast)$. Then the $c$-vectors $c_j^\ast$ of $T^\ast$ are obtained from the $c$-vectors $c_j$ of $T$ by the following recipe.}

(1) \emph{$c_k^\ast=-c_k$}

(2) \emph{$c_j^\ast=c_j+c_k$ if $\ell_j$ connects $p_k$ to its rightmost child or $q_k$ to its leftmost parent.}

(3) \emph{Otherwise, $c_j^\ast=c_j$.}

\noindent\emph{In particular, $T^\ast$ is uniquely determined by $T$ and the edge $\ell=(p_k,q_k)$. We write $T^\ast=\mu_k(T)$.}
\smallskip

{\bf Proof:} In $T$, $p_k$ is a child of $q_k$. When the slope of $\ell_k$ changes from positive to negative, $q_k$ gains an extra parent and $p_k$ loses one parent. To correct this, we move the leftmost parent of $q_k$ over to $p_k$. If this edge is $\ell_j$ then $c_j$ changes to $c_j+c_k$. Similarly, $p_k$ gains one extra child and $q_k$ loses one child. This is corrected by sliding the rightmost child of $p_k$ over to $q_k$. Again the corresponding $c$-vector changes from $c_j$ to $c_j+c_k$. Now each node has the correct number of children and parents. 

{\bf Claim:} This sliding process does not introduce self-intersections of the tree. So, the new tree $T'$ is embedded in the plane.

{\bf Proof:} Suppose, e.g., that sliding the rightmost child of $p_k$, call it $t$, over to $q_k$ produces a self-intersection of the tree. This means that, in the original tree $T$, the triangle with vertices $p_k,q_k,t$ which already has two edges in $T$ meets $T$ at another point $r$ on the third side $tq_k$. By symmetry we may assume that $t$ is to the left of $q_k$. Since $t$ is below $p_k$ which is below $q_k$, the vertical ray going up from $r$ meets the edge $\ell_k$ at some point $s$. By the Claim in the proof of Proposition 1, the unique path in $T$ going from $r$ to $s$ is always increasing which implies that this path goes through $p_k$ but not through $t$ contradicting the assumption that $t$ is the rightmost child of $p_k$. (If $p_k$ has two children, we cannot go from $r$ to $p_k$ through its left child without crossing the wall descending from $p_k$.) The argument for sliding parents is analogous.

By the Uniqueness Lemma, the new tree $T'$ equals $T^\ast$. This proves Proposition 2.

\bigskip

\centerline{\bf Representations of quivers}\smallskip

Associated to the $\e$ vector $\e=(\e_1,\cdots,\e_n)$ we define the \emph{quiver} $Q_\e$ of type $A_{n-1}$ to be the quiver with $n-1$ vertices $1,2,\cdots,n-1$ with each of the $n-2$ pairs of consecutive vertices connected by an arrow:
\[
	Q_{(-, +,-,-,+)}: \bullet \xleftarrow{\ +\ }\bullet\xrightarrow{\ -\ } \bullet\xrightarrow{\ -\ } \bullet
\]
The orientation of the $i$-th arrow (connecting vertex $i$ to vertex $i+1$) should point left if $\e_{i+1}=1$ and right if $\e_{i+1}=-1$. Note that the first and last signs $\e_1,\e_n$ in $\e$ are ignored.

The \emph{Euler matrix} $E_Q$ of any quiver $Q$ without oriented cycles or multiple edges is defined to be the $n-1\times n-1$ square matrix with diagonal entries equal to 1 and all other entries equal to 0 except for the entries $e_{ij}=-1$ when there is an arrow $i\to j$ in $Q$. For $Q=Q_\e$, the Euler matrix will be denoted $E_\e$. For example,
\[
	E_{(-, +,-,-,+)}=\small\mat{
	1&0&0&0\\
	-1&1&-1&0\\
	0&0&1&-1\\
	0&0&0&1
	}
\]

A \emph{representation} of any quiver $Q$ over any field $K$ is given by assigning a finite dimensional $K$-vector space $V_i$ to each vertex $i$ and a $K$-linear map $V_\alpha:V_i\to V_j$ to each arrow $\alpha:i\to j$ in $Q$. An important example is the representation $M_{pq}$ of $Q_\e$ for any $1\le p<q\le n$. This is given by $M_i=(M_{pq})_i=K$ for all $p\le i<q$ with $M_\alpha:M_i\to M_j$ being the identity map on $K$ if $p\le i,j<q$ and $\alpha:i\to j$ is an arrow in $Q_\e$. For any representation $V$, a \emph{subrepresentation} $W$ of $V$ is given by taking a vector subspace $W_i$ of $V_i$ at each vertex $i$ with the property that $V_\alpha(W_i)\subseteq W_j$ for every arrow $\alpha:i\to j$. We need the following lemma.\smallskip

{\bf Lemma:} (Subroot Lemma) 

(a) \emph{If $p<r<q$ then $M_{pr}$ is a subrepresentation of $M_{pq}$ if and only if $\e_r=1$.}

(b) \emph{If $p<r<q$ then $M_{rq}$ is a subrepresentation of $M_{pq}$ if and only if $\e_r=-1$.}

(c) \emph{If $p<a<b<q$ then $M_{ab}$ is a subrepresentation of $M_{pq}$ if and only if $\e_a=-1$ and $\e_b=1$.
}

{\bf Proof:} (a) If $\e_r=1$ then there is an arrow $r-1\xleftarrow\alpha r$ in the quiver $Q_\e$. In the representation $M=M_{pq}$ we have $K$ at both of these points and $M_\alpha=id_K$. So, a subrepresentation cannot be equal to $K$ at $r$ and $0$ at $r-1$. Therefore, $M_{rq}$ is not a subrepresentation of $M_{pq}$ in that case. But $M_{pr}$ would be a subrepresentation of $M_{pq}$. 

(b) If $\e_r=-1$ the arrow goes the other way $\alpha:r-1\to r$. So, $M_{pr}$ cannot be a subrepresentation of $M_{pq}$ but $M_{rq}$ is a subrepresentation of $M_{pq}$. 

(c) In order for $M_{ab}$ to be a subrepresentation of $M_{pq}$ the arrows must be pointing inward at each end: $a-1\to a$ and $b\to b-1$. This is equivalent to $\e_a=-1$ and $\e_b=1$.
\smallskip

Given any representation $V$, the \emph{dimension vector} $\underline{\dim}\,V\in \mathbb N^{n-1}$ of $V$ is the column vector whose $i$-th coordinate is $\dim_KV_i$. In the case of $Q_\e$, the dimension vectors of the representation $M_{pq}$ is equal to
\[
	\underline\dim\,M_{pq}=\beta_{pq}=\sum_{i=p}^{q-1}e_i=\gamma_p-\gamma_q,\qquad 1\le p<q\le n
\]
where $e_i$ is the $i$-th standard unit vector and $\gamma_k=e_k+e_{k+1}+\cdots+e_{n-1}=\beta_{kn}$ if $k<n$ and $\gamma_n=0$. The vectors $\beta_{pq}$ are called the \emph{positive roots} of the quiver $Q_\e$. ($-\beta_{pq}$ are the \emph{negative roots}.) \emph{Note that $c$-vectors of mixed cobinary trees are, by definition, always positive or negative roots.} Also, for any two representations $V,W$ of any quiver $Q$ without oriented cycles we have the well-known formula:
\[
	\dim_K\Hom(V,W)-\dim_K\Ext^1(V,W)=\underline{\dim}\,V^t\, E_Q\,\underline{\dim}\,W.
\]
\bigskip

\centerline{\bf Exchange matrix}\smallskip

If $T$ is a MCT with $c$-matrix $C(T)=C$, let $X_\e=E_\e-E_\e^t$. Then the \emph{exchange matrix} of $T$ is
\begin{equation}\label{eq: def of exchange matrix}
	\tilde B(T)=\mat{C^tX_\e C\\C}
\end{equation}
This formula is a theorem of Nakanishi and Zelevinsky [NZ] which we are using as a definition. The skew symmetric matrix $B=C^tX_\e C$ is called the \emph{principal part} of $\tilde B$. We will show that the matrix $\tilde B(T)$ mutates according to the formula of Fomin and Zelevinsky [FZ]. First, we need to review their definitions.

Suppose that $B, C$ are $n\times n$ integer matrices where $B$ is \emph{skew-symmetrizable} in the sense that its transpose $B^t$ is equal to $-DBD^{-1}$ for some diagonal matrix $D$ with positive integer diagonal entries. The (initial) exchange matrix is defined to be the $2n\times n$ matrix $\tilde B=\mat{B\\C}$.

{\bf Definition} [FZ] For $k=1,\cdots,n$ the \emph{mutation} $\mu_k(\tilde B)$ of $\tilde B=(b_{ij})$ in the \emph{$k$-direction} is defined to be the matrix $\tilde B^\ast=(b'_{ij})$ given as follows.

(1) {$b_{ij}'=-b_{ij}$ if $i=k$ or $j=k$
}

(2) {$b_{ij}'=b_{ij}+ b_{ik}|b_{k j}|$ if $i,j\neq k$ and $b_{ik}$ and $b_{k j}$ have the same sign.
}

(3) {$b_{ij}'=b_{ij}$ otherwise.
}

The composition of a sequence of mutations $\mu_{k_m}\cdots\mu_{k_1}$ will be called a \emph{multiple mutation}. The \emph{$c$-vectors} corresponding to a skew symmetrizable matrix $B$ are defined to be those vectors which occur as columns of the bottom half $C'$ of matrices $\tilde B'$ obtained from $\tilde B=\mat{B\\I_n}$ by a multiple mutation.\smallskip

Nakanishi and Zelevinski proved in [NZ] that in many case, such as when $B$ is skew-symmetric, any iterated mutation of $\tilde B=\mat{B\\I_n}$ has the form $\tilde B'=\mat{C^tBC\\C}$ which is \eqref{eq: def of exchange matrix} with $X_\e=B$. This is the sense in which we are taking their theorem as a definition.\smallskip

{\bf Theorem 2:} \emph{
Let $T^\ast=\mu_k(T)$. Then the exchange matrix $\tilde B'=\tilde B(T^\ast)$ is equal to the mutation of $\tilde B=\tilde B(T)$ in the $k$-direction. I.e.: 
}
\[
	\tilde B(\mu_k(T))=\mu_k(\tilde B(T)).
\]

{\bf Remark 2:} Let $T_0$ be the unique MCT with $\pi(T_0)=\{id\}$ (the second tree in Figure 1 and the fourth tree in Figure 4). Then $B(T_0)=\mat{X_\e\\I_{n-1}}$. So, Theorem 2 implies that our geometric notion of a $c$-vector for a MCT $T$ (by definition the columns of the lower half of $\tilde B(T)$) agrees with the Fomin-Zelvinski notion of $c$-vector for the skew symmetrix matrix $X_\e$ (the columns of the lower half of a multiple mutation of $\tilde B(T_0)$). To draw this conclusion we also need Remark~1 above which implies that every MCT $T$ is a multiple mutation of $T_0$.
\smallskip

{\bf Proof:} We reduce the number of cases to check by noting that the statement is invariant under vertical and horizontal symmetry of $T$: Under vertical symmetry, $C,X_\e,B$ become $-C,-X_\e,-B$. Under rotation by $\pi$, the order of the rows and columns of $C, X_\e$ and $B$ are reversed. The matrix $X_\e=E_\e-E_\e^t$ is skew-symmetric. The only nonzero entries are $\e_i$ along the superdiagonal and $-\e_i$ along the subdiagonal for $1<i<n$:
\[
	X_\e=E_\e-E_\e^t=\small\mat{0& \e_2&0&0&\cdots\\
	-\e_2&0&\e_3&0&\\
	0&-\e_3&0&\e_4&\ddots\\
	0&0&-\e_4&0&\ddots\\
	\vdots&&\ddots&\ddots&\ddots
	}
\]
From this we have the following easy calculation:
\[
	\gamma_p^t X_\e \gamma_q=\begin{cases} \e_{q} & \text{if } p<q\\
   -\e_{p} & \text{if } p>q\\
   0 & \text{if } p=q
    \end{cases}
\]
which also implies that $\beta_{pr}^t X_\e\beta_{qr}=\gamma_p^t X_\e \gamma_q$ if $p,q<r$.

We will now calculate $b_{ij}=c_i^t X_\e c_j$ for $i,j<n$ in all possible cases.

{\bf Case 1:} If $c_i=\pm\beta_{pq}$ and $c_j=\pm\beta_{p'q'}$ are $c$-vectors for $T$ with $p,q,p',q'$ distinct integers then $c_i^tX_\e c_j=0$.

{\bf Proof:} Up to symmetry, there are three cases

(1) $p<q<p'<q'$. Then $c_i^tX_\e c_j=c_i^tX_\e \gamma_{p'}-c_i^tX_\e \gamma_{q'}=0-0=0$.

(2) $p<p'<q<q'$. Then $c_i^tX_\e \gamma_{q'}=\e_{q'}-\e_{q'}=0$. So, 
\[
	c_i^tX_\e c_j=c_i^tX_\e \gamma_{p'}=\gamma_p^tX_\e \gamma_{p'}-\gamma_q^tX_\e \gamma_{p'}=\e_{p'}+\e_{q}=0
\]
since $\e_{p'}$, $\e_{q}$ must have opposite sign: If $\ell_i$ lies above $\ell_j$ then $\e_{p'}=1$ and $\e_{q}=-1$ by definition of mixed cobinary tree. If $\ell_i$ lies below $\ell_j$ then $\e_{p'}=-1$ and $\e_{q}=1$.

%
\begin{figure}[htbp]
\begin{center}
%
{
\setlength{\unitlength}{1cm}
{\mbox{
\begin{picture}(6,2.5)
      \thicklines
\qbezier(0,.5)(1.5,1.25)(3,2)
\spot{0,.5}
\spot{4,0}
\spot{3,2}
\spot{2,.5}
\qbezier(2,.5)(3,.25)(4,0)
\thinlines
\swline{2,.5}
\veeup{3,2}
\qbezier(2,.5)(3.5,1)(5,1.5)
\put(-0.2,.7){$p$}
\put(1,1.3){$\ell_i$}
\put(1.8,.7){$p'$}
\put(2.7,-.1){$\ell_j$}
\put(3.2,1.8){$q$}
\put(4.1,.1){$q'$}
\end{picture}}
}}
\caption{If $p<p'<q<q'$ and $\ell_i$ is above $\ell_j$ then $\e_{p'}=1$ and $\e_q=-1$.}
\label{figureC}
\end{center}
\vspace{-12pt}
\end{figure}
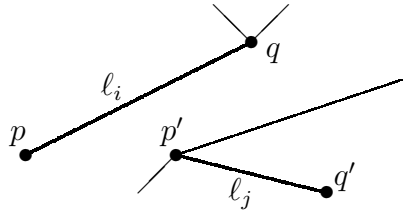
(3) $p<p'<q'<q$. Then $c_i^tX_\e c_j=\gamma_p^tX_\e c_j=\e_{p'}-\e_{q'}=0$ since $\e_{p'}=\e_{q'}$: If $\ell_i$ lies above $\ell_j$ they are both positive, if $\ell_i$ lies below $\ell_j$ they are both negative.

{\bf Case 2:} Now suppose that the internal edges $\ell_j,\ell_k$ of $T$ have the same right endpoint or the same left endpoint. Then $|b_{kj}|=1$ and $\sgn (b_{kj})=\sgn (c_k)$.

{\bf Proof:} The edges $\ell_j,\ell_k$ must have opposite sign. By symmetry we may assume $\ell_j$ is negative and $\ell_k$ is positive. Then $c_k=\beta_{pr}$ and $c_j=-\beta_{qr}$ for some $p,q<r$. The formula for $b_{kj}$ is
\[
	b_{kj}=c_k^t X_\e c_j=- \gamma_p^t X_\e \gamma_q=\begin{cases} -\e_q & \text{if } p<q\\
    \e_p & \text{if } q<p
    \end{cases}
\]
In either case, $b_{kj}=1$ as claimed. The case when $\ell_j,\ell_k$ share a left endpoint is easier.

{\bf Case 3:} Suppose that the left endpoint of $\ell_j$, say $t_q$, is equal to the right endpoint of $\ell_k$. Then
\begin{equation}\label{eq for bkj}
	b_{kj}=\e_q\, \sgn (c_j)\,\sgn (c_k)
\end{equation}
In particular, $b_{kj}$ and $c_k$ have the same sign if and only if $\sgn   (c_j)=\e_q$.

{\bf Proof:} If $c_j,c_k$ have the same sign then
\[
	b_{kj}=c_k^t X_\e c_j=(\gamma_p-\gamma_q)^tX_\e(\gamma_q-\gamma_r)=\gamma_p^t X_\e \gamma_q=\e_q
\]
If $c_j,c_k$ have opposite sign then $b_{kj}=-\gamma_p^t X_\e \gamma_q=-\e_q$. In both cases \eqref{eq for bkj} holds.

{\bf Case 4:} If the right endpoint of $\ell_j$, say $t_q$, is equal to the left endpoint of $\ell_k$ then, by reversing $j,k$ we get: $b_{kj}=-b_{jk}=-\e_q\sgn (c_j)\,\sgn (c_k)$. 

This exhausts all possible configurations of $\ell_j$ and $\ell_k$ and we saw that $|b_{kj}|\le 1$ for all $k,j$. If $\ell_k$ goes from $t_p$ to $t_q$ (so that $c_k=\pm \beta_{pq}$) then we have the following.

{\bf Lemma.} \emph{For each $j<n$, $b_{kj}$ is nonzero with the same sign as $c_k$ if and only if one of the following holds. Each condition has two versions depending on the sign of $c_k$.}

(a) (For $c_k$ positive) \emph{$\ell_j$ connects $t_q$ to the leftmost parent of $t_q$.} (For $c_k$ negative) \emph{$\ell_j$ connects $t_p$ to the rightmost parent of $t_p$.}

(b) (For $c_k$ positive) \emph{$\ell_j$ connects $t_p$ to the rightmost child of $t_p$.} (For $c_k$ negative) \emph{$\ell_j$ connects $t_q$ to the leftmost child of $t_q$.}\smallskip

If we compare this with Proposition 2, we see that mutation of $T$ in the $k$-direction changes the $c$-matrix $C$ of $T$ by column operations as follows.

(1) Add Column $k$ to Column $j$ whenever $b_{kj}$ is $\pm1$ with the same sign as $c_k$.

(2) Change the sign of Column $k$.

This shows that the lower half of the matrix $\tilde B(T)$ which is $C$ mutates according to the given formula due to Fomin and Zelevinsky. The principal part, which is $B=C^t X_\e C$, mutates by the same column operations and the corresponding row operations. Since $B$ is skew-symmetric, these row operations are given as follows.

(3) Add Row $k$ to Row $j$ whenever $b_{jk}$ is $\pm1$ with the opposite sign as $c_k$.

(4) Change the sign of Row $k$.

Whenever $b_{ik}=1$ and $b_{kj}=-1$, the operations (1) and (3) cancel each other and $b_{ij}$ is unchanged. If $b_{ik}=-1$ and $b_{kj}=1$ then nothing is done to either Column $i$ or Row $j$, so $b_{ij}$ is again unchanged. If $b_{ik}$ and $b_{kj}$ are equal then either Row operation (3) will add $b_{ik}$ to $b_{ij}$ or Column operation (1) will add $b_{kj}$ to $b_{ij}$. The result is the same. Thus the tree mutation $\mu_k$ changes $\tilde B(T)$ by matrix mutation $\mu_k$ as required, proving the theorem.
\bigskip

\centerline{\bf Definition of a cluster}\smallskip

To explain the statement that mixed cobinary trees represent clusters of representations of the quiver $Q_\e$, we first need to define clusters. We will do this in an ultimately elementary way by first reviewing the theoretical definition, then reducing the definition of a cluster to a matrix equation. In the next section we will show that MCT's produce solutions of this equation.

Before giving definitions we recall that the indecomposable representations of $Q_\e$ are $M_{pq}$ which are uniquely determined by their dimension vectors which are the positive roots $\beta_{pq}=\gamma_p-\gamma_q$ where $1\le p<q\le n$. Every representation of $Q_\e$ can be expressed uniquely up to permutation of summands as a direct sum of these representations. The dimension vectors of the indecomposable projective representations $P_i$ are the \emph{projective roots} $\pi_i$. These are given as the rows of the matrix $E_\e^{-1}$. (The matrix $E_\e$ is $I_{n-1}$ minus a nilpotent matrix $X$ with nonnegative entries. Thus $E_\e^{-1}=I_{n-1}+X+X^2+\cdots+X^{n-2}$ has nonnegative entries.) We define an \emph{almost positive root} to be either a positive root $\beta_{pq}$ or a negative projective root $-\pi_i$. There are $\binom{n}2+n$ almost positive roots for $Q_\e$.

The \emph{cluster category} introduced in [BMRRT] has as (indecomposable) objects the indecomposable representations of $Q_\e$ plus the shifted projective representations $P_i[1]$ with dimension vectors $-\pi_i$. A \emph{cluster} in the cluster category can be defined to be a set of $n-1$ distinct objects $M_i$ from this collection of objects so that $\Ext^1(M_i,M_j)=0$ for all $i,j$. If $\underline{\dim} M_i=v_i$, this is equivalent to the condition that
\begin{equation}\label{numerical characterization of cluster}
v_i^t E_\e v_j\ge0
\end{equation}
since $v_i^t E_\e v_j=\dim_K\Hom(M_i,M_j)-\dim_K\Ext^1(M_i,M_j)$ and $\Hom(M_i,M_j),\Ext^1(M_i,M_j)$ cannot both be nonzero (since $Q_\e$ is Dynkin). More precisely, a set of $n-1$ distinct almost positive roots $v_i$ form a cluster if and only if \eqref{numerical characterization of cluster} holds for all $i,j$ so that $v_j$ is a positive root.

We define a \emph{cluster matrix} for $Q_\e$ to be an $(n-1)\times (n-1)$ matrix $V$ whose columns are distinct almost positive roots $v_i$ so that all of the entries of $V^tE_\e W$ are nonnegative where $W$ is the submatrix of $V$ whose columns are those $v_i$ which are positive. (This definition work whenever the quiver is a Dynkin diagram.) The \emph{initial cluster matrix} is defined to be $(E_\e^t)^{-1}$, the matrix whose columns are the dimension vectors $\pi_i$ of the projective representations $P_i$ of $Q_\e$. This is a cluster matrix since $(E_\e)^{-1}E_\e (E_\e^t)^{-1}=(E_\e^t)^{-1}$. We note that the columns of a cluster matrix give an ordered cluster. But, the cluster itself is unordered.

One important fact that we need is that the cluster matrix $V$ is always invertible.\smallskip

{\bf Lemma:} \emph{The objects in a cluster can be ordered in such a way that $V^t E_\e V$ is unipotent, i.e., upper triangular with 1's on the diagonal. In particular $V$ is invertible as an integer matrix.}

{\bf Proof:} In the case at hand, $Q_\e$ is a Dynkin diagram. So, we can arrange the objects of the cluster according to their position in the Auslander-Reiten (AR) quiver. Since the AR quiver has the property that there are no nonzero morphism going from right to left, there will be no nonzero morphisms $M_j\to M_i$ for $j>i$. So, $v_j E_\e v_i=\dim_K \Hom(M_j,M_i)=0$. Furthermore, $\underline\dim M^t E_\e\underline\dim M=1$ for all indecomposable objects $M$. So, the statement hold.

In general, this lemma follows from Theorem 2.4 in [S92] which says that the modules in a cluster can be arranged so that $\Hom(M_j,M_i)=0$ for $j>i$. The shifted projective modules $P_i[1]$ in any cluster, not considered in [S92], should be taken after the modules and so that $\Hom(P_j,P_i)=0$ for $j>i$. Then the matrix $V^t E_\e V$ will be upper triangular. Also, all objects in a cluster are real Schur representations ([K], [S92], [IOTW]) which implies that $\dim \Hom(M_i,M_i)=v_i^t E_\e v_i=1$.  So, the lemma holds for any quiver $Q$ without oriented cycles.


\bigskip

\centerline{\bf Clusters and weights of semi-invariants}\smallskip

Hugh Thomas explained the following formula \eqref{HughThomasEquation} to the first author. (See [ST].) 
\begin{equation}\label{HughThomasEquation}
	V^tE_Q C=I_{n-1}
\end{equation}
This formula gives the relation between a cluster matrix $V$ and the corresponding $c$-matrix $C$. In order for this to make sense as a theorem, we would need to explain this correspondence. Instead, we will use this formula as the definition of the $c$-vectors corresponding to a cluster. We define $C:=(E_Q)^{-1}(V^t)^{-1}$ to be the \emph{classical $c$-matrix} corresponding to a cluster matrix $V$. By the lemma above this is an integer matrix. Since the transpose of a permutation matrix is equal to its inverse any permutation of the columns of $V$ will permute the columns of $C$ in the same way. Therefore, this formula assigns a \emph{classical $c$-vector} $c_i'$ to every object in every (unordered) cluster. However, the same object may have several associated $c$-vectors if it lies in several clusters.

In the case of the quiver $Q_\e$, we will show that this classical $c$-matrix is equal to the $c$-matrix given by the geometry of mixed cobinary trees as explained in this paper. To do this, we will use virtual semi-invariants associated to a cluster by [IOTW]. We will use the following elementary matrix version of the results of [IOTW].

Suppose that $Q$ is a quiver whose underlying graph is a Dynkin diagram with $n-1$ vertices, for example, $Q=Q_\e$. Recall that for every positive root $\beta$ of $Q$ there is an indecomposable module $M_\beta$ with dimension vector $\beta$. Also recall that $\beta'$ is a subroot of $\beta$ if $M_\beta$ contains a subrepresentation isomorphic to $M_\beta$. (See the Subroot Lemma.) Define $D(\beta)$ to be set of all $v\in  \mathbb R^{n-1}$ satisfying the following two conditions.

(a) $v^t E_Q \beta=0$

(b) (\emph{stability condition}) $v^t E_Q \beta'\le 0$ for all subroots $\beta'$ of $\beta$.

In [IOTW], $D(\beta)$ is called the \emph{support} of the \emph{virtual semi-invariant} $\sigma_\beta$ with \emph{weight} $E_Q\beta$.

{\bf Example 1:} Let $Q_\e$ be the quiver $1\to 2\to 3$ and let $\beta=(0,1,1)^t$. Then $E_\e\beta=(-1,0,1)^t$. So, for any $v^t=(x,y,z)$, Condition (a) states $x=z$. The corresponding semi-invariant $\sigma_\beta$ is the determinant of the map $M_1\to M_3$ for any representation $M$ of $Q_\e$ which is defined when $\dim M_1=\dim M_3$ or, equivalently, $\underline\dim M^t E_\e\beta=0$. This is a semi-invariant since it depends on a choice of bases for the vectors spaces $M_1,M_3$.

The only subroot of $\beta$ is $\beta'=(0,0,1)^t$ with $E_\e\beta'=(0,-1,1)^t$. So, stability condition (b) is: $y\ge z$. This corresponds to the fact that $\dim M_2\ge \dim M_3$ is a necessary condition in order for the semi-invariant $\sigma_\beta(M)=\det(M_1\to M_3)$ to be nonzero.

The ``determinantal'' semi-invariant $\sigma_\beta(M)$ is given in general as follows ([S91]). Let $0\to P_1\to_f P_0\to M\to 0$ be an exact sequence where $P_0,P_1$ are projective representations. Then $\sigma_\beta(M)$ is the determinant of the induced map:
\[
	f^\ast:\Hom(P_1,M_\beta)\to \Hom(P_0,M_\beta).
\]
By [IOTW], there exists a morphism $f:P_1\to P_0$ between projective representations of $Q$ so that $f^\ast$ is an isomorphism if and only if $v=\dim P_0-\dim P_1$ satisfies (a) and (b) above. Such a morphism $f$ is called a \emph{virtual representation} and the extension of $\sigma_\beta$ to such $f$ is called a \emph{virtual semi-invariant}.

Given any cluster matrix $V$, let $S(V)\subseteq \mathbb R^{n-1}$ be the set of all $\sum a_iv_i$ where $a_i\ge0$. Then the Virtual Stability Theorem of [IOTW] can be stated as follows.
\smallskip

{\bf Theorem 3:}
\emph{Let $V=(v_i)$ be a cluster matrix for the Dynkin quiver $Q$. Then there exist unique positive roots $\beta_1,\cdots,\beta_{n-1}$ of $Q$ so that $v_j\in D(\beta_i)$ for $i\neq j$. Furthermore, there is a triangulation of the unit sphere $S^{n-2}$ in $\mathbb R^{n-1}$ whose $n-2$ simplices are $S^{n-2}\cap S(V)$ for all cluster matrices $V$ and so that the $n-3$ skeleton of the triangulation is equal to the union of all $D(\beta)\cap S^{n-2}$ for all positive roots $\beta$ of $Q$.} 
\smallskip


{\bf Remark 3:} The statement $v_j\in D(\beta_i)$ has two parts. The first part is (a) $v_j^t E_Q \beta_i=0$ for $j\neq i$. Comparing this with Equation \eqref{HughThomasEquation} we see that $\beta_i$ must be an integer multiple of the classical $c$-vector $c_i'$ for each $i$. This integer must be $\pm 1$. So, $c_i'=\pm \beta_i$ for each $i$.\smallskip



Using the stability condition (b) $v_j^t E_\e\beta'\le0$ for $\beta'\subset \beta_i$, we will show that there is a bijection between clusters for the quiver $Q_\e$ and mixed cobinary trees $T$ having the property that the classical $c$-vectors of the cluster are equal to the $c$-vectors of the corresponding MCT. This bijection uses the following equation for the closed cone $\overline R(T)\subset \mathbb R^n$ associated to $T$. Let $F:\mathbb R^n\to \mathbb R^{n-1}$ be given by
\[
	F(x)=(x_2-x_1,x_3-x_2,\cdots,x_n-x_{n-1}).
\]
Then $\overline R(T)$ is the set of all $x\in\mathbb R^n$ so that $c_i^t F(x)\ge 0$ for every $c$-vector $c_i$ associated to $T$. Equivalently:
\[
	\overline R(T)=F^{-1}\{
	y\in \mathbb R^{n-1}\,|\, y^t c_i\ge 0\text{ for all $c$-vectors $c_i$ of $T$}
	\}
\]

\smallskip

{\bf Theorem 4:} \emph{Given any cluster $\{M_1,\cdots,M_{n-1}\}$ in the cluster category for $Q_\e$ with associated cluster matrix $V$, there is a unique mixed cobinary tree $T$ with the property that}
\[
	\overline R(T)=F^{-1}(E_\e^t S(V)).
\]
Since $E_\e^t S(V)$ is the set of all $y\in \mathbb R^{n-1}$ so that $y^tc_i'\ge 0$ for all classical $c$-vectors $c_i'$ associated to $V$, we conclude that $c_i'=c_i$ (up to renumbering of the $c$-vectors):\smallskip

{\bf Corollary 2:} \emph{Up to permutation, the $c$-vectors $c_i$ associated to a mixed cobinary tree $T$ are equal to the classical $c$-vectors $c_i'$ associated to the corresponding cluster $(v_i)$. Equivalently, 
\[
v_i^tE_\e c_j=\delta_{ij}.
\]}

{\bf Remark 4:} This says that the mixed cobinary trees label the faces of the cluster complex for the quiver $Q_\e$ and therefore also label the vertices of the dual object which is the generalized associahedron for $Q_\e$. (See [Ste].) 

{\bf Example 2:} Let $n=5$ with $\e=(-1,1,-1,1,1)$. Then $Q_\e= \bullet\leftarrow\bullet\rightarrow\bullet\leftarrow\bullet$ and $E_\e$ and its inverse are given by:
\[
	E_\e=\small\mat{1 &0&0&0\\
	-1&1&-1&0\\
	0&0&1&0\\
	0&0&-1&1}\qquad
	E_\e^{-1}=\mat{
	1&0&0&0\\
	1&1&1&0\\
	0&0&1&0\\
	0&0&1&1
	}
\]
The rows of the matrix $E_\e^{-1}$ are the projective roots. In particular $\pi_4=(0,0,1,1)^t$. Consider the cluster whose matrix is $V$ which, together with $V^t E_\e$, are:
\[
	V=\small\mat{
	1&1&0&0\\
	1&1&1&0\\
	1&0&1&-1\\
	0&0&0&-1
	}\qquad
	V^tE_\e=\mat{
	0&1&0&0\\
	0&1&-1&0\\
	-1&1&0&0\\
	0&0&0&-1
	}
\]
$V$ is a cluster matrix since its columns are positive roots and negative projective roots, e.g. $v_4=-\pi_4$, and all entries of $V^tE_\e W$ are nonnegative where $W$ is the $4\times 3$ matrix consisting of the first three columns of $V$.

Theorem 4 tells us to take the set of all nonnegative linear combinations of the rows of $V^t E_\e$ and apply $F^{-1}$. Up to addition of a scalar multiple of the vector $(1,1,1,1,1)$, $F^{-1}$ of the rows of $V^t E_\e$ gives:
\[
	\small\mat{
	0&0&1&1&1\\
	0&0&1&0&0\\
	1&0&1&1&1\\
	1&1&1&1&0
	}
\]
The sum of the rows is $(2,1,4,3,2)$. Applying the algorithm in the proof of the Uniqueness Lemma (for either $21543$ or $31542$), we get the following tree:
%
\begin{figure}[htbp]
\begin{center}
%
{
\setlength{\unitlength}{1cm}
{\mbox{
\begin{picture}(5,2.5)
      \thicklines
\veedown{1,.5}
\veedown{4,1}
\veeup{2,2}
\nwline{0,1}
\qbezier(0,1)(.5,.75)(1,.5)
\qbezier(0,1)(1.5,1.25)(3,1.5)
\qbezier(2,2)(3,1.5)(4,1)
\put(-0.3,.7){$t_1$}
\put(.8,1.4){$\ell_1$}
\put(3.2,.9){$\ell_4$}
\put(.2,.4){$\ell_3$}
\put(2.4,1.9){$\ell_2$}
\put(3.1,1.6){$t_4$}
\put(1.1,.6){$t_2$}
\put(1.7,1.7){$t_3$}
\put(4.1,1.1){$t_5$}
\end{picture}}
}}
\label{figureD}
\end{center}
\vspace{-12pt}
\end{figure}

The four edges of this tree give the inequalities $x_2\le x_1\le x_4\le x_3$ and $x_5\le x_4$. By definition $\overline R(T)$ is the set of all $x\in\mathbb R^n$ satisfying these inequalities. The statement of Theorem 4 is that $\overline R(T)$ is equal to the set of nonnegative linear combinations of the rows of the above $4\times 5$ matrix plus multiples of $(1,1,1,1,1)$. To see that this holds in our example note that, if we cut the edge $\ell_1$, the lower end is connected only to $t_1,t_2$ and the upper end is connected to $t_3,t_4,t_5$ giving the vector $(0,0,1,1,1)$ which is Row 1. Thus, increasing the slope of $\ell_1$ will add a multiple of Row 1 to the vector in $\overline R(T)$. Cutting the other three edges gives the other three rows of the matrix.

The matrix of $c$-vectors of $T$ is:
\[
	C=\small\mat{
	1&0&-1&0\\
	1&0&0&0\\
	1&-1&0&0\\
	0&0&0&-1
	}
\]
which is equal to $(V^tE_\e)^{-1}$ as stated. Thus, Theorem 4 holds in  Example 2.

We list the four cluster objects in Example 2, their associated $c$-vectors, corresponding semi-invariants (applied to a representation $M$) and stability conditions in terms of $\underline\dim M=(w,x,y,z)^t$. We use the shorthand $[w\,x\,y\,z]:=(w,x,y,z)^t$ to save space in this table.\smallskip
\[
\text{Example 2: one cluster in } Q_{-+-++}\]\[
\begin{array}{|c|c|c|c|c|}
\hline
\text{cluster} &\text{$c$-vectors} &\text{weights}& \text{semi-invariants}  &\text{stability conditions}\\
\text{objects} & c_i=c_i'&\text{of } \sigma_i&&\underline\dim M E_\e\beta_i= 0\\
v_i &  =\pm\beta_i & E_\e \beta_i& \sigma_i=\sigma_{\beta_i}(M)& \beta ,\quad \underline\dim M E_\e\beta\le 0
\\
\hline
&&&& w+y=x+z
\\
{[1\,1\,1\,0]} & [1\,1\,1\,0] & [1,\text{-}1,1,\text{-}1] & \det(M_2\oplus M_4\to M_1\oplus M_3) & [1\,0\,0\,0],\quad  w\le x
\\
&&&&  [0\,0\,1\,0],\quad  y\le x+z
\\
\hline
[1\,1\,0\,0] & -[0\,0\,1\,0] & [0,\text{-}1,1,\text{-}1] &\det(M_2\oplus M_4\to  M_3) & y=x+z
\\
\hline
[0\,1\,1\,0] &-[1\,0\,0\,0] & [1,\text{-}1,0,0] &\det(M_2\to M_1)&w=x
\\
\hline
-[0\,0\,1\,1] &-[0\,0\,0\,1] & [0\,0\,0\,1] & \det(0\to M_4)&z=0
\\
\hline
\end{array}
\]\smallskip
For $i=1$, the two subroots of $\beta_1= [1\,1\,1\,0] $ are $ [1\,0\,0\,0] $ and $ [0\,0\,1\,0] $. The corresponding stability conditions $w\le x$ and $y\le x+z$ hold for any nonnegative linear combination of $v_2,v_3,v_4$. The condition $w+y=x+z$ also holds on such linear combinations. Therefore, $D(\beta_1)$ contains all nonnegative linear combinations of $v_2,v_3,v_4$. The other roots $\beta_2,\beta_3,\beta_4$ are simple and thus have no subroots. For $i=4$ we use the convention that $\det (0\to 0)=1$. Therefore, the semi-invariant $\sigma_{\beta_4}(M)$ is nonzero exactly when $z=\dim M_4=0$.
\smallskip

{\bf Proof of Theorem 4:} The theorem is equivalent to the statement that $E_\e^t S(V)=F\overline R(T)=F(\overline R(T)\cap H)$ where $H$ is the set of all $x\in\mathbb R^n$ with $\sum x_i=0$. We will show that, for every cluster matrix $V$, there exists an MCT $T$ so that $E_\e^t S(V)\subseteq F\overline R(T)$. Since the sets $S(V)$ cover $\mathbb R^{n-1}$ and $E_\e$ is an invertible matrix, this implies that each set $F\overline R(T)$ is a union of sets $E_\e^t S(V)$. So, we get a surjective mapping from the set of clusters in the cluster category of the quiver $Q_\e$ to the set of MCT's $T$. By Theorem 1 and the well-known properties of clusters of type $A_{n-1}$, both sets have a Catalan number $C_n$ of elements. Therefore, any epimorphism is a bijection and we conclude that $E_\e^t S(V)=F\overline R(T)$, proving the theorem.

Both $E_\e^t S(V)$ and $F\overline R(T)=F(\overline R(T)\cap H)$ are cones on $n-2$ simplices. Their unions are cones on partitions of the unit $n-2$ sphere in $\mathbb R^{n-1}$. The regions $E_\e^t S(V)$ give a triangulation of the sphere. Showing that the top dimensional simplices of the $E_\e^tS(V)$ cluster triangulation are subsets of the regions $F(\overline R(T)\cap H)$ is equivalent to the statement that the codimension-1 simplices of the latter (the image under $F$ of the union of the walls of the $\overline R(T)$'s) is a subset of the $n-3$ skeleton of the cluster triangulation. By Theorem 3, this is the union of the supports $D(\beta)$. Thus, all we need to do is to show that each wall of each $\overline R(T)\cap H$ is contained in $E_\e^t D(\beta)$ for some positive root $\beta$.


Thus, we start with an MCT $T$ with $c$-vectors $c_i=\delta_i \beta_{p_i q_i}$ where $\delta_i=sgn(\ell_i)$. Consider a general point $x=(x_1, \cdots,x_n)$ on the $j$-th face of $\overline R(T)$ given by the equation $x_{p_j}=x_{q_j}$ and the inequalities $\delta_i (x_{q_i}-x_{p_i})\ge0$ for $i\neq j$. In terms of $y=F(x)$, these conditions are 

($*$) $y^tc_j=0$ and $y^tc_i\ge0$ for $i\neq j$. 

\noindent By Theorem 3, it suffices to show that $D(\beta_{p_jq_j})$ contains all $z$ where $y=E_\e^t z$. By definition of $D(\beta)$ this is equivalent to the following two conditions where we write $p=p_j,q=q_j$.

(a) $z^tE_\e \beta_{pq}=0$ or, equivalently, $y^t\beta_{pq}=0$. 

(b) \emph{Stability condition:} $y^t \beta_{ab}\le 0$ for every subroot $\beta_{ab}\subseteq \beta_{pq}$.

\noindent But, (a) holds by assumption. To prove (b) we use the Subroot Lemma which says that $\beta_{ab}\subseteq \beta_{pq}$ when $p\le a<b\le q$ and both of the following conditions are satisfied.

(1) Either $a=p$ or $\e_a=-1$.

(2) Either $b=q$ or $\e_b=1$.

\noindent Condition (1) implies that $x_a\ge x_p=x_q$ since $x$ is in the closure of $R(T)$ and $x_a>\min(x_p,x_q)$ for all $x\in R(T)$. Similarly, Condition (2) implies that $x_b\le x_p=x_q$. Therefore $y^t\beta_{ab}=y_a+\cdots+y_{b-1}=x_b-x_a\le0$. So, the stability conditions (b) are satisfied at all $y=F(x)$. So each wall of $F\overline R(T)$ for each $T$ lies in the $n-3$ skeleton of the cluster triangulation which is what we needed to prove.
\bigskip

\centerline{\bf Summary}\smallskip

Mixed cobinary trees are trees $T$ which are embedded in the plane with $n$ internal nodes each of which has either two children and one parent (a $\Lambda$-node) or one child and two parents (a V-node). There is also a left-right order of the internal nodes and a ``wall separation'' condition as illustrated in color in Figure 5. The sign sequence $\e\in\{\pm1\}^n$ gives the left-to-right order of the $\Lambda$ (positive) and V (negative) nodes. A root of type $A_{n-1}$ is a vector in $\{0,1,-1\}^{n-1}$ all of whose nonzero entries are equal and in consecutive positions. Each internal edge of the tree $T$ gives a positive or negative root of the root system of type $A_{n-1}$ whose endpoints are given by the endpoints of the edge and whose sign is given by the slope of the edge. The $c$-matrix $C=C(T)$ is the square matrix whose columns are the $c$-vectors of $T$ in any order. The region $R(T)\subset \mathbb R^n$ is the set of all $x\in \mathbb R^n$ so that $y^tc_i>0$ for each $c$-vector $c_i$ where $y=F(x)\in\mathbb R^{n-1}$ is the vector with coordinates $y_i=x_{i+1}-x_i$. This is equal to the set of all points $x\in \mathbb R^n$ so that the $n$ points in the plane $(i,x_i)$ form allowable positions for the internal nodes of $T$. We defined mutation $T^\ast=\mu_k(T)$ if $\overline R(T)\cap \overline R(T^\ast)$ is the face of $\overline R(T)$ given by $y^tc_k=0$.

Given $\e$, there is an associated quiver $Q_\e$ of type $A_{n-1}$ using only the signs $\e_2,\cdots,\e_{n-1}$ to indicate the direction of the arrows between the $n-1$ points in the quiver ($\e=1$ being left and $\e=-1$ being right). $E_\e$ is the Euler matrix of $Q_\e$. It has 1's on the diagonal and $ij$ entry equal to $-1$ if there is an arrow $i\to j$ in the quiver. Then $B(T)=C^t(E_\e-E_\e^t)C$ is the principal part of the exchange matrix $\tilde B(T)=\mat{B(T)\\C(T)}$. We showed that this matrix transforms according to the formula of Fomin and Zelevinsky, i.e., $\tilde B(\mu_k(T))=\mu_k\tilde B(T)$.

The projective roots $\pi_i$ of $Q_\e$ were defined to be the rows of the matrix $E_\e^{-1}$. A cluster was defined to be a set of $n-1$ distinct roots which are either positive roots or negative projective roots so that $V^tE_\e W$ has nonnegative entries where $V$ is the square matrix whose columns $v_i$ are the elements of the cluster and $W$ is the submatrix of $V$ whose columns are the positive roots of the cluster. The classical $c$-matrix was defined by the equation $C=(V^tE_\e)^{-1}$. In the final theorem, the formula $E_\e^tS(V)=F\overline R(T)$ was shown to give a bijection between clusters and trees so that the classical $c$-vectors of $V$ are equal to the $c$-vectors of the corresponding tree $T$.\bigskip

\centerline{\bf Acknowledgements}\smallskip

The relation between $c$-vectors and semi-invariants was explained to the first author (in different language) by both Hugh Thomas and Salvatore Stella. The formula relating internal edges of mixed binary trees with weights of semi-invariants (and thus with $c$-vectors) is a special case of joint work with Gordana Todorov, Jerzy Weyman and Kent Orr. The definition of mixed cobinary trees arises naturally from the Stability Theorem for virtual semi-invariants. Many years ago, Robert Penner explained to the first author how triangulations of surfaces can be represented by locally planar graphs called \emph{fat graphs} [P]. These are also called \emph{ribbon graphs}. Mixed cobinary trees are examples of fat graphs. The first author is partially supported by National Security Agency Grant \# H98230-13-1-0247.
\bigskip

\centerline{\bf References}\smallskip

[BMRRT] Aslak~Bakke Buan, Robert~J. Marsh, Idun Reiten, and Gordana Todorov,
  \emph{Tilting theory and cluster combinatorics}, Adv. Math. \textbf{204}
  (2006), no.~2, 572--618.

[FZ] Sergey Fomin and Andrei Zelevinsky, \emph{Cluster algebras. {IV}.
  {C}oefficients}, Compos. Math. \textbf{143} (2007), no.~1, 112--164.

[IOTW] Kiyoshi Igusa, Kent Orr, Gordana Todorov, and Jerzy Weyman, \emph{Cluster
  complexes via semi-invariants}, Compos. Math. \textbf{145} (2009), no.~4,
  1001--1034.
  
[K]  V. G. Kac, \emph{Infinite root systems, representations of graphs and invariant theory. II}, J. Algebra 78 (1982), no. 1, 141Ð162.
  
[TZ] Tomoki Nakanishi, Andrei Zelevinsky, \emph{On tropical dualities in cluster algebras}, arXiv: 1101.3736.

[P]
Robert~C. Penner, \emph{The decorated {T}eichm\"uller space of punctured
  surfaces}, Comm. Math. Phys. \textbf{113} (1987), no.~2, 299--339.

[S91] Aidan Schofield, \emph{Semi-invariants of quivers}, J. London Math. Soc. (2) \textbf{43} (1991), no.~3, 385--395.
 
[S92] Aidan Schofield, \emph{General Representations of quivers}, Proc. London Math. Soc. (3) \textbf{65} (1992), 46-64.

[ST] David Speyer and Hugh Thomas, \emph{Acyclic cluster algebras revisited}, arXiv:1203.0277.


[Ste] Salvatore Stella, \emph{Polyhedral models for generalized associahedra via Coxeter elements}, arXiv:111.1657.


\end{document}